\documentclass[12pt,a4paper]{amsart}
\usepackage{amsmath, amssymb, amsfonts,enumerate}
\usepackage[all]{xy}
\usepackage{amscd}

\newtheorem{theorem}{Theorem}[section]
\newtheorem{lemma}[theorem]{Lemma}
\newtheorem{corollary}[theorem]{Corollary}
\newtheorem{proposition}[theorem]{Proposition}

\newtheorem*{question*}{Question}

\theoremstyle{definition}

\theoremstyle{remark}

\numberwithin{equation}{section}

\newcommand {\K}{\mathbb{K}} %% a field
\newcommand {\N}{\mathbb{N}} %% positive integers
\newcommand {\Z}{\mathbb{Z}} %% integers
\DeclareMathOperator{\Ker}{Ker}

\DeclareMathOperator{\Id}{Id}

\begin{document}
\title[On the reversibility and the closed image property]{On the reversibility and the closed image property of linear cellular automata}
\date{\today}

\author{Tullio Ceccherini-Silberstein}
\address{Dipartimento di Ingegneria, Universit\`a del Sannio, C.so
Garibaldi 107, 82100 Benevento, Italy}
\email{tceccher@mat.uniroma1.it}
\author{Michel Coornaert}
\address{Institut de Recherche Math\'ematique Avanc\'ee,
UMR 7501,                                             Universit\'e  de Strasbourg et CNRS,
                                                 7 rue Ren\'e-Descartes,
                                               67000 Strasbourg, France}
\email{coornaert@math.unistra.fr}
\date{\today}
\begin{abstract}
When $G$ is an arbitrary  group and $V$ is a finite-dimensional vector space,
it is known that every bijective linear cellular automaton
$\tau \colon V^G \to V^G$ is reversible
and that the image of every linear cellular automaton $\tau \colon V^G \to V^G$ is closed in $V^G$ for the prodiscrete topology. 
In this paper, we present a new proof of these two results which is based on the Mittag-Leffler lemma for projective sequences of sets. 
We also show that if $G$ is a non-periodic group and $V$ is an infinite-dimensional vector space, then there exist a  linear cellular automaton $\tau_1 \colon V^G \to V^G$ which is bijective but not reversible and a  linear cellular automaton $\tau_2 \colon V^G \to V^G$ whose image is not closed in $V^G$ for the  prodiscrete topology.
 \end{abstract}

\maketitle

% SECTION 1
\section{Introduction}
\label{sec:introduction}

 Let $G$ be a group and let $A$ be a set.
The set $A^G$, which consists of all maps $x  \colon G \to A$ is called the set of \emph{configurations} over the group $G$ and the \emph{alphabet} $A$.
There is a natural left action of the group $G$ on $A^G$ defined by $gx(h) = x(g^{-1}h)$ for all 
$g,h \in G$ and $x \in A^G$.
This action is called the $G$-\emph{shift } on $A^G$ and   is the fundamental object of study in the branch of mathematics known as \emph{symbolic dynamics}.
Although classical symbolic dynamics often restricts to the case when the alphabet set $A$ is finite,
  it is clear from recent developments in the theory, such as the ones contained in the influential work of M. Gromov   \cite{gromov-esav},  \cite{gromov-tids},
that the study of shift systems with infinite alphabet $A$ also deserves attention, especially when 
$A$ is equipped with some additional (linear, algebraic, symplectic, etc.) structure.
 \par
A \emph{cellular automaton} over the group $G$ and the alphabet $A$ is a map 
$\tau \colon A^G \to A^G$ satisfying the following property:
there exist a finite subset $M \subset G$ 
and a map $\mu \colon  A^M \to A$ such that 
\begin{equation} 
\label{e;local-property}
\tau(x)(g) = \mu(\pi_M(g^{-1}x))  \quad  \text{for all } x \in A^G \text{ and } g \in G,
\end{equation}
where $\pi_M \colon A^G \to A^M$ denotes the restriction map.
Such a set $M$ is then called a \emph{memory set} for $\tau$ and $\mu$ is called the \emph{local defining map} for  $\tau$ associated with $M$. It follows from this definition that every cellular automaton $\tau \colon A^G \to A^G$ commutes with the shift action, i.e., it satisfies 
\begin{equation}
\label{e:G-equivariance}
\tau(gx) = g\tau(x) \ \ \mbox{$g \in G$ and $x \in A^G$.}
\end{equation}
In other words, $\tau$ is $G$-equivariant.
\par
A  cellular automaton $\tau \colon A^G \to A^G$ is said to be \emph{reversible} if $\tau$ is bijective and  the inverse map $\tau^{-1} \colon A^G \to A^G$
is also a cellular automaton.
\par
We equip $A^G$ with its \emph{prodiscrete} topology, that is, with the product topology obtained by taking the discrete topology on each factor $A$ of $A^G$.
This turns out $A^G$ into a totally disconnected Hausdorff topological space.
Every cellular automaton $\tau \colon A^G \to A^G$ is continuous with respect to the prodiscrete topology. Conversely, when the alphabet $A$ is finite, it follows from the Curtis-Hedlund theorem \cite{hedlund}  that every continuous map $f \colon A^G \to A^G$ which commutes with the shift action is a cellular automaton.
\par
When $A$ is finite, the space $A^G$ is compact by the Tychonoff product theorem and one immediately deduces from the Curtis-Hedlund theorem 
 that every bijective cellular automaton 
$\tau \colon A^G \to A^G$ is reversible. 
  \par
Suppose now that $V$ is a vector space over a field $\K$.
A cellular automaton $\tau \colon V^G \to V^G$ over the alphabet $V$  is said to be \emph{linear} if 
$\tau$ is $\K$-linear with respect to the natural $\K$-vector space structure on $V^G$.
In this setting we have the following:

\begin{theorem}
\label{t:reversibilite}
Let $G$ be a group and let $V$ be a finite-dimensional vector space over a field $\K$.
Then every bijective linear cellular automaton $\tau \colon V^G \to V^G$ is reversible.
\end{theorem} 

This result was proved in \cite{israel} under the assumption that the group $G$ is countable and then extended to any group $G$ in \cite{induction}.  
In this paper, we use the Mittag-Leffler lemma for projective sequences of sets
to derive a new proof of Theorem \ref{t:reversibilite}.
We then provide examples of bijective linear cellular automata which are not reversible. More precisely, we shall prove the following (recall that a group is said to be \emph{periodic} if all its elements have finite order):

\begin{theorem}
\label{t:lca-not-reversible} 
Let $G$ be a non-periodic group and let $\K$ be a field.
Let $V$ be an infinite-dimensional vector space over $\K$. 
Then there exists a bijective linear cellular automaton $\tau \colon V^G \to V^G$ which is not reversible.   
\end{theorem}

If $A$ is an infinite set, it is always possible to find a vector space with the same cardinality as $A$.
One can take for example the $\K$-vector space based on $A$ for an arbitrary finite field $\K$. Therefore, as an immediate consequence of Theorem \ref{t:lca-not-reversible},
we get: 

\begin{corollary}
\label{c:ca-not-reversible}
Let $G$ be a non-periodic group and let $A$ be an infinite set.
Then there exists a bijective cellular automaton $\tau \colon A^G \to A^G$ which is not reversible.\qed
\end{corollary}

Given a set $X$ and a topological space $Y$, one says that a map $f \colon X \to Y$  has the 
\emph{closed image property} if the set $f(X)$ is closed in $Y$.
The closed image property is often used to establish surjectivity results.
Indeed, to prove that a map $f \colon X \to Y$ with the closed image property is surjective,
it suffices to show that $f(X)$ is dense in $Y$. 
 When the alphabet $A$ is finite, every cellular automaton 
$\tau \colon A^G \to A^G$ has the closed image property
(with respect to the prodiscrete topology) 
by compactness of $A^G$. 
In the linear setting one has:

\begin{theorem}
\label{t:closed-image}
Let $G$ be a group and let $V$ be a finite-dimensional vector space over a field $\K$.
Then every linear cellular automaton $\tau \colon V^G \to V^G$ has the closed image property
with respect to the prodiscrete topology on $V^G$. 
\end{theorem} 

This was proved in \cite{garden} when the group $G$ is countable, and then extended to any group $G$ in \cite{induction} (see also \cite[Section 4.D]{gromov-esav} for more general results).
As for  the reversibility result mentioned above (cf. Theorem \ref{t:reversibilite}),  we shall present here a new proof of Theorem \ref{t:closed-image} based on the Mittag-Leffler lemma for projective sequences of sets. Moreover, we shall also establish the following:

\begin{theorem}
\label{t:lca-not-closed-image}
Let $G$ be a non-periodic group and let $\K$ be a field.
Let $V$ be an infinite-dimensional vector space over $\K$. 
Then there exists a linear cellular automaton $\tau'\colon V^G \to V^G$ such that 
 $\tau'(V^G)$ is not closed in $V^G$ with respect to the prodiscrete topology.
\end{theorem}

As above, this gives us:

\begin{corollary}
\label{c:not-closed} 
Let $G$ be a non-periodic group and let $A$ be an infinite set. 
Then there exists a cellular automaton $\tau'\colon A^G \to A^G$ such that $\tau'(A^G)$ is not closed in $A^G$ with respect to the prodiscrete topology.\qed
\end{corollary}

The remainder of the paper is organized as follows. In Section \ref{sec:induction},
for the sake of completeness and the convenience of the reader, we briefly recall from  \cite{induction} the definitions of induction and restriction of cellular automata and list some of their  properties.
 Section \ref{sec:ML} is devoted to the Mittag-Leffler lemma for projective sequences of sets  . This     set-theoretic version of the Mittag-Leffler lemma may be easily deduced 
 from Theorem 1 in \cite[TG II. Section 5]{bourbaki-top-gen} 
 (see also \cite[Section I.3]{grothendieck-ega-3}) but we present here a self-contained   proof for the convenience of the reader.  
The proofs of  Theorem \ref{t:reversibilite} and Theorem \ref{t:closed-image} which are based on the 
Mittag-Leffler lemma
are given in Section \ref{sec:proofs-finite-dim}.
Each proof is divided into two steps.
We first establish the result in the case when $G$ is countable by means of the Mittag-Leffler lemma and then extend it to the general case  by applying the first step to the cellular automaton obtained by restriction to the  subgroup generated by a memory set.
In Section \ref{sec:proofs}, we give the proofs of Theorem \ref{t:lca-not-reversible} and Theorem \ref{t:lca-not-closed-image}. These are also divided into two steps : we first treat the case $G = \Z$ and then we use the technique of induction   to extend them to any non-periodic group.
In the final section, we discuss the question whether the results from Section \ref{sec:proofs} extend to periodic groups which are not locally finite. 
\par
This research was initiated in November 2008 during our participation to the Conference ``ESI Workshop on Structural Probability'' held at the Erwin Schr\"{o}dinger Institute in Vienna.  We would like to express our gratitude to the organizers for their hospitality and the nicest stimulating atmosphere at the ESI.

% SECTION 2
\section{Induction and restriction of cellular automata}
\label{sec:induction}

In this section we recall the notions of induction and restriction for cellular automata
(cf. \cite{induction}).
\par
Let $G$ be a group, $A$ a set, and $H$ a subgroup of $G$.
\par
Suppose that a cellular automaton $\tau \colon A^G \to A^G$ admits a memory set $M$ such 
that $M \subset H$. Let $\mu \colon A^M \to A$ denote the associated local defining map.
Then the map $\tau_H \colon A^H \to A^H$ defined by
$$
\tau_H(y)(h) = \mu(\pi_M(h^{-1}y))
\quad \text{ for all } y \in A^H, h \in H,
$$
is a cellular automaton over the group $H$ and the alphabet $A$ with memory set $M$ and local defining map $\mu$.
One says that $\tau_H$ is the cellular automaton obtained by \emph{restriction} of  
$\tau$ to $H$.
\par
Conversely, let $\sigma \colon A^H \to A^H$  be a cellular automaton with memory set $M \subset H$ and local defining map
$\mu \colon A^M \to A$. Then the map $\sigma^G \colon A^G \to A^G$ defined by
$$
\sigma^G(x)(g) = \mu(\pi_M(g^{-1}x))
\quad \text{ for all } x \in A^G, g \in G,
$$
is a cellular automaton over the group $G$ and the alphabet $A$ with memory set $M$ and local defining map $\mu$.
One says that $\sigma^G$ is the cellular automaton obtained by \emph{induction} of $\sigma$ to $G$.
\par
It immediately follows from their definitions that induction and restriction are operations one inverse to the other in the sense that one has $(\tau_H)^G = \tau$ 
and $(\sigma^G)_H = \sigma$ for every cellular automaton  $\tau \colon A^G \to A^G$ over $G$ admitting a memory set contained in $H$ and every cellular automaton $\sigma \colon A^H \to A^H$ over $H$.
We shall use the following results (see \cite[Theorem 1.2]{induction} for proofs):

\begin{theorem}
\label{t:induction}
Let $G$ be a group, $A$ a set, and $H$ a subgroup of $G$. Let $\tau \colon A^G \to A^G$ be a cellular automaton admitting
a memory set contained in $H$. Then the following holds.
\begin{enumerate}[{\rm (i)}]
\item $\tau$ is bijective if and only if $\tau_H$ is bijective;
\item $\tau$ is reversible if and only if $\tau_H$ is reversible; 
\item $\tau(A^G)$ is closed in $A^G$ if and only if $\tau_H(A^H)$ is closed in $A^H$ 
(for the prodiscrete topology on $A^H$);
\item when $A$ is a vector space, $\tau$ is linear if and only if $\tau_H$ is linear.
\end{enumerate}
\end{theorem}

% SECTION 3
\section{The Mittag-Leffler lemma} 
\label{sec:ML}

In this section, we give the proof of the version of the Mittag-Leffler lemma 
that we shall use in the next section in order to establish
  Theorem \ref{t:reversibilite} and Theorem \ref{t:closed-image}.
Let us first recall a few facts about projective limits of projective sequences in the category of sets.
\par
Denote by $\N$ the set of nonnegative integers.
A \emph{projective sequence} of sets is a sequence $(X_n)_{n \in \N}$ of sets equipped with maps $f_{nm} \colon X_m \to X_n$, defined for all $n,m \in \N$ with $ m \geq n$,
satisfying the following conditions:

\begin{enumerate}[(PS-1)]
\item
$f_{n n}$ is the identity map on $X_n$ for all $n \in \N$;
\item
$f_{n k} = f_{n m} \circ f_{m k}$ for all $n,m,k \in \N$ such that $k \geq m \geq n$.
\end{enumerate}
We denote such a projective sequence by $(X_n,f_{n m})$ or simply by $(X_n)$.
The \emph{projective limit} $\varprojlim X_n$ of the projective sequence $(X_n,f_{n m})$ is the subset of $\prod_{n \in \N} X_n$ consisting of the sequences $(x_n)_{n \in \N}$ satisfying 
$x_n = f_{n m}(x_m)$ for all $n,m \in \N$ such that $m \geq n$.
\par
One says that the projective sequence $(X_n)$ satisfies the \emph{Mittag-Leffler condition} if the following holds:
\begin{itemize}
\item[(ML)]{for each $n \in \N$, there exists $m \in \N$ with $m \geq n$ such that $f_{n k}(X_k) = f_{n m}(X_m)$ for all $k \geq m$.}
\end{itemize}

\begin{lemma}[Mittag-Leffler]
\label{l;ML}
Let  $(X_n,f_{n m})$ be a projective sequence of nonempty sets
satisfying the Mittag-Leffler condition. Then its projective limit $X = \varprojlim X_n$ is not empty.
\end{lemma}

\begin{proof}
We first observe that given an arbitrary projective sequence of sets $(X_n,f_{n m})$,
then Property (PS-2) implies that, for each $n \in \N$, the sequence of sets $f_{nm}(X_m)$, $m \geq n$, is non-increasing. Let us set $X_n' = \bigcap_{m \geq n} f_{nm}(X_m)$ (this is called the set of \emph{universal elements} in $X_n$, cf. \cite{grothendieck-ega-3}).
The map $f_{n m}$ clearly induces by restriction a map $g_{n m} \colon X_m' \to X_n'$ for all $m \geq n$. Then $(X_n',g_{n m})$ is a projective sequence 
having the same projective limit as the projective sequence $(X_n,f_{n m})$.
\par
Suppose now that all the sets $X_n$ are nonempty and that the projective sequence $(X_n,f_{n m})$ satisfies the Mittag-Leffler condition.
This means that, for each $n \in \N$, there is an integer $m \geq n$ such that $f_{n k}(X_k) = f_{n m}(X_m)$ for all $k \geq m$.
This implies $X_n' = f_{n m}(X_m)$ so that, in particular, the set $X_n'$ is not empty.
We claim that the map $g_{n, n + 1} \colon X_{n + 1}' \to X_n'$ is surjective for every $n \in \N$. To see this, let $n \in \N$ and $x_n' \in X_n'$. By the Mittag-Leffler condition, we can find an integer $p \geq n + 1$ such that
$f_{n k}(X_k) = f_{n p}(X_p)$ and $f_{n + 1,  k}(X_k) = f_{n + 1,p}(X_p)$ for all $k \geq p$. 
It follows that $X_n' = f_{n p}(X_p)$ and $X_{n + 1}' = f_{n + 1, p}(X_p)$. Consequently, we can find $x_p \in X_p$ such that $x_n' = f_{n p}(x_p)$.
Setting $x_{n + 1}' = f_{n + 1,p}(x_p)$, we have  $x_{n + 1}' \in X_{n + 1}'$ and 
$$
g_{n, n + 1}(x_{n +1}') = f_{n, n + 1}(x_{n +1}') = f_{n, n + 1} \circ f_{n+1,p} (x_p) = f_{n p}(x_p) = x_n'.
$$ 
This proves our claim that $g_{n, n + 1}$ is onto.
As the sets $X_n'$ are nonempty,   we can now construct by induction a sequence $(x_n')_{n \in \N}$ such that $x_n' = g_{n, n + 1}(x_{n +1}') $ for all $n \in \N$.
This sequence is in the projective limit $\varprojlim X_n' = \varprojlim X_n$. This shows that $\varprojlim X_n$ is not empty.    
\end{proof}

% SECTION 4
\section[Reversibility and the closed image property]{Reversibility and the closed image property in the finite-dimensional case}
\label{sec:proofs-finite-dim}

This section contains the proofs of   Theorem \ref{t:reversibilite} and Theorem \ref{t:closed-image}.
based on the Mittag-Leffler lemma for projective sequence of sets.

 \begin{proof}[Proof of Theorem \ref{t:reversibilite}]
 Let $\tau \colon V^G \to V^G$ be a bijective linear cellular automaton.
 We have to show that $\tau$ is reversible. 
We split the proof into two steps.
\par 
Suppose first that the group $G$ is countable. 
Since $\tau$ is linear and $G$-equivariant (cf. \eqref{e:G-equivariance}), the inverse map $\tau^{-1} \colon V^G \to V^G$ is linear and $G$-equivariant as well. 
Let us show that the following local property is satisfied
by $\tau^{-1}$: there exists a finite subset $N \subset G$ such that
\begin{itemize}
\item[($\ast$)] for $y \in V^G$, the element $\tau^{-1}(y)(1_G)$ only depends on
the restriction of $y$ to $N$.
\end{itemize}
This will show that $\tau$ is reversible.
Indeed, if ($\ast$) holds for some finite subset $N \subset G$, then there exists a (unique) map $\nu \colon V^{N} \to V$ satisfying
$$
\tau^{-1}(y)(1_G) = \nu(\pi_{N}(y)).
$$
From the $G$-equivariance of $\tau^{-1}$ we then deduce
$$
\tau^{-1}(y)(g) = g^{-1}\tau^{-1}(y)(1_G) = \tau^{-1}(g^{-1}y)(1_G) =  \nu(\pi_N(g^{-1}y)).
$$
for all $y \in V^G$, which implies that $\tau^{-1}$ is the cellular automaton with memory set $N$ and local defining map $\nu$.
\par
Let us assume by contradiction that there exists no finite subset $N \subset G$  satisfying condition ($\ast$). 
Let $M$ be a memory set for $\tau$ such that $1_G \in M$.
Since $G$ is countable, we can find a sequence $(A_n)_{n \in \N}$ of finite subsets of $G$ such that $G = \bigcup_{n \in \N} A_n$,
$M \subset A_0$ and $A_n \subset A_{n + 1}$ for all $n \in \N$. 
Let $B_n = \{g \in G: gM \subset A_n\}$. Note that $G = \bigcup_{n \in \N} B_n$,
$1_G \in B_0$,  and $B_n \subset B_{n + 1}$
for all $n \in \N$.
\par
Since there exists no finite subset $N \subset G$ satisfying condition ($\ast$), we can find, for each $n \in \N$, two configurations $y_n', y_n'' \in V^G$ such that $y_n'\vert_{B_n} = y_n''\vert_{B_n}$ and ${\tau}^{-1}(y_n')(1_G)
\neq {\tau}^{-1}(y_n'')(1_G)$. By linearity of ${\tau}^{-1}$, the configuration $y_n = y_n' - y_n'' \in V^G$
satisfies 
\begin{equation}
\label{e:noyeau-non-vide}
y_n\vert_{B_n} = 0 \quad \text{ and } \quad \tau^{-1}(y_n)(1_G) \neq 0.
\end{equation}
  
It follows from \eqref{e;local-property} that if $x$ and $x'$ are elements in $V^G$ such that $x$ and $x'$ coincide on $A_n$  then the configurations $\tau(x)$ and $\tau(x')$ coincide on $B_n$. Therefore, given $x_n \in V^{A_n}$ and denoting by $\widetilde{x_n} \in V^G$ a configuration extending $x_n$,
the element
$$
u_n = \tau(\widetilde{x_n})\vert_{B_n} \in V^{B_n}
$$ 
does not depend on the particular choice of the extension $\widetilde{x}_n$ of $x_n$. Thus we can define a map
$\tau_n \colon V^{A_n} \to V^{B_n}$ by setting $\tau_n(x_n) = u_n$. It is clear that $\tau_n$ is $\K$-linear.
\par
Consider, for each $n \in \N$, the subset $X_n \subset V^{A_n}$ 
consisting of all $x_n \in V^{A_n}$ such that $x_n \in \Ker(\tau_n)$ and $x_n(1_G) \neq 0$. 
Note that $X_n$ is not empty since $(\tau^{-1}(y_n))\vert_{A_n} \in X_n$ by 
\eqref{e:noyeau-non-vide}. Now observe that, for $m \geq n$, the restriction 
map  $\rho_{n m} \colon V^{A_m} \to V^{A_n}$ is $\K$-linear and induces a   map 
$f_{nm} \colon X_m \to X_n$.
Indeed, if $u \in X_m$, then we have $u\vert_{A_n} \in X_n$ since   
$\tau_n(u\vert_{A_n}) = (\tau_m(u))\vert_{B_n} = 0$ and $(u\vert_{A_n})(1_G) = u(1_G) \neq 0$.
Conditions (PS-1) and (PS-2) are trivially satisfied so that
$(X_n,f_{nm})$ is a projective sequence of nonempty sets. 
Let us show that $(X_n,f_{nm})$
also satisfies the Mittag-Leffler condition (ML).
Consider, for all $m \geq n$, the set $f_{nm}(X_m) \subset X_n \subset V^{A_n}$.
By definition, we have that $f_{nm}(X_m) = \rho_{n m}(\Ker(\tau_m)) \cap X_n$.
Observe now that, if $n \leq m \leq m'$, then
$\rho_{n m'}(\Ker(\tau_{m'})) \subset \rho_{n m}(\Ker(\tau_m))$.
 Therefore, if we fix $n$, the sequence 
$\rho_{n m}(\Ker(\tau_m))$, where $m=n, n+1,\ldots$, is a non-increasing sequence of vector subspaces of $V^{A_n}$. 
 As the vector space $ V^{A_n}$ is finite-dimensional, this sequence stabilizes, and it follows that, for each 
 $n \in \N$, there exists an integer $m \geq n$ such that $f_{nk}(X_k) = f_{nm}(X_m)$ if $k \geq m$.  This shows that condition (ML) is satisfied.
We then deduce from Lemma \ref{l;ML} that the projective limit $X =\varprojlim X_n$ is not empty. 
Let $(z_n) \in X  $. Then there exists a unique $z \in V^G $ such that $\pi_{A_n}(z) = z_n$ for all 
$n \in \N$. But $\tau(z) = 0$ since $\pi_{B_n}(\tau(z)) = \tau_n(z_n) = 0$ for all $n $
and $z(1_G) = z_0(1_G) \not= 0$. 
 This contradicts the injectivity of $\tau$.
\par
This shows, that there exists a finite subset $N \subset G$ satisfying ($\ast$) and therefore that $\tau$ is reversible.
\par
We now drop the countability assumption on $G$ and prove the theorem in the general case. Choose a memory set $M \subset G$ for $\tau$ and denote by $H$ the   subgroup of $G$ generated by  $M$.
Observe that $H$ is countable since $M$ is finite.
  By assertions (i) and (iv) of Theorem \ref{t:induction},
   the restriction cellular automaton $\tau_H \colon V^H \to V^H$ is linear and bijective. It then follows from the previous step that $\tau_H$ is reversible (that is, the inverse map $(\tau_H)^{-1} \colon V^H \to V^H$ is a cellular automaton). 
   By applying assertion (ii) of Theorem \ref{t:induction}, we conclude that $\tau$ is also reversible.
\end{proof}

\begin{proof}[Proof of Theorem \ref{t:closed-image}]
Let $\tau \colon V^G \to V^G$ be a linear cellular automaton.
We have to show that $\tau$ has the closed image property with respect to the prodiscrete topology on $V^G$.
We split the proof into two steps as in the preceding proof.
\par 
Suppose first that the group $G$ is countable.
Choose, as in the first step of the preceding proof,   a sequence $(A_n)_{n \in \N}$ of finite subsets of $G$ such that $G = \bigcup_{n \in \N} A_n$,
$M \subset A_0$ and $A_n \subset A_{n + 1}$,   and consider, for each $n \in \N$, 
the $\K$-linear map $\tau_n \colon V^{A_n} \to V^{B_n}$, where
 $B_n = \{g \in G: gM \subset A_n\}$
 and $\tau_n$ is defined by $\tau_n(x_n) = (\tau(\widetilde{x_n}))\vert_{B_n}$ for all $x_n \in V^{A_n}$ and $\widetilde{x_n} \in V^G$ extending $x_n$.
\par
 Let now $y \in V^G$ and suppose that $y$ is in the closure of $\tau(V^G)$.
Then, for all $n \in \N$, there exists $z_n \in V^G$ such that 
\begin{equation}
\label{e;y-z-n-B-n}
\pi_{B_n}(y)=  \pi_{B_n}(\tau(z_n)).
\end{equation}
Consider, for each $n \in \N$, the affine
subspace $X_n \subset V^{A_n}$ defined by $X_n = \tau_n^{-1}(\pi_{B_n}(y))$.
We have $X_n \not= \varnothing$ for all $n$ by \eqref{e;y-z-n-B-n}. 
For $m \geq n$, the restriction map $V^{A_m} \to V^{A_n}$ induces an affine map
$f_{nm} \colon X_m \to X_n$. Conditions (PS-1) and (PS-2) are trivially satisfied so that
$(X_n,f_{nm})$ is a projective sequence. We claim that $(X_n,f_{nm})$
also satisfies the Mittag-Leffler condition (ML). Indeed, consider, for all $m \geq n$, the
affine subspace $f_{nm}(X_m) \subset X_n$. We have $f_{nm'}(X_{m'})  \subset f_{nm}(X_m)$ for all $n \leq m \leq m'$ since $f_{nm'} = f_{nm} \circ f_{mm'}$.
As the sequence  $f_{nm}(X_m)$ ($m = n,n+1,\dots$) is a
non-increasing sequence of finite-dimensional affine subspaces, it
stabilizes, i.e., for each $n \in \N$ there exists an integer $m \geq n$ such that
$f_{nk}(X_k) = f_{nm}(X_m)$ if $k \geq m$. Thus, condition (ML) is satisfied.
It follows from Lemma \ref{l;ML} that the projective limit $\varprojlim X_n$ is nonempty.
Choose an element $(x_n)_{n \in \N} \in \varprojlim X_n$. We have that
$x_{n + 1}$ coincides with $x_n$ on $A_n$ and that $x_n \in V^{A_n}$ for all $n \in \N$. 
As $G = \cup_{n \in \N} A_n$, we deduce that there exists a (unique) configuration 
$x \in V^G$ such that $x\vert_{A_n} = x_n$ for all $n$. 
We have $\tau(x)\vert_{B_n}= \tau_n(x_n) = y_n = y\vert_{B_n}$ for all $n$. 
Since $G = \cup_{n \in \N} B_n$, this shows that $\tau(x) = y$.
This completes the proof in the case that $G$ is countable.
\par
Let us treat now the case of an arbitrary (possibly uncountable) group $G$.
As in the second step of the preceding proof,
choose a memory set $M \subset G$ for $\tau$
and consider the countable  subgroup of $G$ generated by $M$.
   By the previous step, we have that the restriction cellular automaton $\tau_H \colon V^H \to V^H$ has the closed image property, that is, $\tau_H(V^H)$
is closed in $V^H$ for the prodiscrete topology.  By applying Theorem \ref{t:induction}.(iii), we deduce that $\tau(V^G)$ is also closed in $V^G$ for the prodiscrete topology.  Thus $\tau$ satisfies the closed image property.
\end{proof}

% SECTION 5
\section{Proofs of Theorem \ref{t:lca-not-reversible} and of  Theorem \ref{t:lca-not-closed-image}}
\label{sec:proofs}

In this section we present the proofs of Theorem \ref{t:lca-not-reversible} and Theorem \ref{t:lca-not-closed-image}. In both cases, we divide the proofs into two steps: we first treat the case $G = \Z$ and then we use the technique of induction  to extend them to any non-periodic group.
\par
Let $\K$ be a field and let $V$ be an infinite-dimensional vector space over $\K$.
\par
Since $V$ is infinite-dimensional, we can find an infinite sequence $(v_i)_{i \geq 1}$ 
of linearly independent vectors in $V$.
Let $E$ denote the vector subspace spanned by all the vectors $v_i$   and let $F$ be a vector subspace of $V$ such that $V = E \oplus F$.
\par
Let us first construct a bijective linear cellular automaton $\sigma \colon V^\Z \to V^\Z$ which is not reversible.
\par
For each $j \geq 1$, denote by $E_j$ the vector subspace spanned by the vectors $v_i$, 
where $(j - 1)j/2 + 1 \leq  i \leq j(j + 1)/2$.
 We have
$E = \bigoplus_{j \geq 1}E_j$ and $\dim_\K (E_j) = j$ for all $j \geq 1$.
\par
For each $j \geq 1$, let $\varphi_j \colon E_j \to E_j$ denote the unique $\K$-linear map such that
 $$
 \varphi_j(v_i) =
 \begin{cases}
 0 &\text{ if } i = (j -1)j/2 + 1, \\
 v_{i -1} & \text{ if } (j - 1)j/2 + 2 \leq i \leq j(j + 1)/2. \\
 \end{cases}
 $$
Observe that $\varphi_j$ is nilpotent of degree $j$.
\par
Consider now the maps $\sigma_j \colon E_j^\Z \to E_j^\Z$ defined by
$$
\sigma_j(x_j)(n) = x_j(n) - \varphi_j(x_j(n + 1)) 
$$
for all $x_j \in E_j^\Z$ and $n \in \Z$.
\par
 We define the map $\sigma \colon V^\Z \to V^\Z$ by
$$
\sigma = \left(\bigoplus_{j \geq 1} \sigma_j\right) \oplus \Id_{F^\Z},
$$
where we use the natural identification $V^\Z = \left(\bigoplus_{j \geq 1} E_j^\Z\right) \oplus F^\Z$ 
and denote by $\Id_X$ the identity map on a set  $X$.

\begin{lemma}
\label{l:sigma}
The map $\sigma \colon V^\Z \to V^\Z$ is a linear cellular automaton which is bijective but not reversible.
\end{lemma}

\begin{proof}
The $\K$-linearity of the maps $\sigma_j$ and hence of $\sigma$  is straightforward from their  definition.
It is also clear that $\sigma$ is a cellular automaton admitting $S = \{0,1\}$ as a memory set. 
The associated local defining map is the map
$\mu \colon V^S \to V$
given by
$$
\mu = \left(\bigoplus_{j \geq 1} \mu_j \right) \oplus \Id_{F^S},
$$
where $\mu_j \colon E_j^S \to E_j^S$ is defined by  $\mu_j(u_0,u_1) = u_0 - \varphi_j(u_1)$ for all $(u_0,u_1) \in E_j^S = E_j \times E_j$.
\par
  Consider the map $\nu_j \colon E_j^\Z \to E_j^\Z$ defined by
\begin{equation}
\label{e:formula-for-nuj}
\nu_j(x_j)(n) = \sum_{k = 0}^{j - 1}  \varphi_j^k(x_j(n + k))
\quad \text{for all } x_j \in E_j^\Z  \text{ and } n \in \Z,
 \end{equation}
 where $\varphi_j^0 = \Id_{E_j}$ and $\varphi_j^{k + 1} = \varphi_j \circ \varphi_j^k$ for all $k \geq 0$.   Using the fact that $\varphi_j$ is nilpotent of degree $j$, one immediately checks that 
 $$
 \sigma_j \circ \nu_j = \nu_j \circ \sigma_j = \Id_{E_j^\Z}.
 $$
Thus  $\sigma_j$ is bijective with inverse map $\nu_j$.
This implies that $\sigma$ is bijective with inverse map
\begin{equation}
\label{e:inverse-of-sigma}
\sigma^{-1} =\left(\bigoplus_{j \geq 1} \nu_j\right) \oplus \Id_{F^\Z}.
\end{equation}
If $\sigma^{-1}$ was a cellular automaton, there would be a finite subset
$M \subset \Z$ such that $\sigma^{-1}(x)(0)$ would only depend of the restriction of $x \in V^\Z$ to $M$.
To see that this is impossible, suppose that there is such an $M$ and choose an integer $j_0 \geq 1$ large enough so that $M \subset (-\infty,j_0 - 2]$.
Consider the configurations $y,z \in V^\Z$ defined by
$y(n) = 0$ for all $n \in \Z$ and 
$$
z(n) = 
\begin{cases}
v_{j_0(j_0 - 1)/2} & \text{ if } n = j_0 - 1,\\
 0 & \text{ if } n \in \Z \setminus \{j_0 - 1\}. \\
\end{cases}
$$
Then $y$ and $z$ coincide on $(-\infty,j_0 - 2]$ and hence on $M$. However, $\sigma^{-1}(y)(0) \not= \sigma^{-1}(z)(0)$ since $\sigma^{-1}(y) = 0$ while, by \eqref{e:inverse-of-sigma} 
and \eqref{e:formula-for-nuj}, we have
$$
\sigma^{-1}(z)(0) = \varphi_{j_0}^{j_0 - 1}(v_{j_0(j_0 + 1)/2}) = v_{(j_0 - 1)j_0/2 + 1} \not= 0.
$$
This shows that $\sigma^{-1}$ is not a cellular automaton. Therefore $\sigma$ is not reversible.
\end{proof}

\begin{proof}[Proof of Theorem \ref{t:lca-not-reversible}]
Since $G$ is non-periodic, we can find an element $g_0 \in G$ with infinite order. 
Denote by $H$ the infinite cyclic subgroup of $G$ generated by $g_0$. 
We identify $\Z$ with $H$ via the map $n\mapsto g_0^n$.
Consider the cellular automaton $\sigma \colon V^H \to V^H$ studied in Lemma \ref{l:sigma}. Then, by virtue of Lemma \ref{l:sigma} and assertions (iv), (i) and (ii) in Theorem  \ref{t:induction}, the  
induced cellular automaton $\tau = \sigma^G \colon V^G \to V^G$ is linear, bijective and not reversible.
\end{proof}

We now construct a linear cellular automaton $\sigma' \colon V^\Z \to V^\Z$ such that $\sigma'(V^\Z)$ is not closed in $V^\Z$ with respect to the prodiscrete topology.
\par
Let $\psi \colon V \to V$ denote the unique $\K$-linear map satisfying $\psi(v_i) = v_{i + 1}$ for all $i \geq 1$, and whose restriction to $F$ is identically $0$.
Define the map $\sigma' \colon V^\Z \to V^\Z$ by
$$
\sigma'(x)(n) = x(n + 1) - \psi(x(n))
$$
for all $x \in V^\Z$ and $n \in \Z$.

\begin{lemma}
\label{l:sigma-2}
The map $\sigma' \colon V^\Z \to V^\Z$ is a linear cellular automaton and $\sigma'(V^\Z)$ is not closed in $V^\Z$.
\end{lemma}

\begin{proof}
The $\K$-linearity of $\sigma'$ immediately follows from the $\K$-linearity of $\psi$.
It is also clear that $\sigma'$ is a cellular automaton admitting $S = \{0,1\}$ as a memory set. 
The associated local defining map is the map
$\mu' \colon V^S \to V$
given by
$\mu'(u_0,u_1) = u_1 - \psi(u_0)$ 
 for all $(u_0,u_1) \in  V^S = V \times V$.
 \par
Let $c \in V^\Z$ denote the constant configuration defined by $c(n) = v_1$  for all $n \in \Z$.
Let us show that $c$ is in the closure of $\sigma'(V^{\Z})$.
By definition of the prodiscrete topology, it suffices to show that for each finite subset $F \subset \Z$, there exists a configuration $x_F \in V^\Z$ such that $c$ and $\sigma'(x_F)$ coincide on $F$. 
To see this, choose an integer $n_0 \in \Z$ such that 
$F \subset [n_0,\infty)$.
Consider the configuration $x_F \in V^\Z$ defined by
$$
x_F(n) =
\begin{cases}
0 & \text{ if } n \leq n_0 - 1 , \\
v_1 + v_2 + \cdots + v_{n - n_0 + 1} & \text{ if } n \geq n_0.  \\
\end{cases}
$$
   Observe that $\sigma'(x_F)(n) = v_1$ for all $n  \geq n_0$,  so that the configurations $\sigma'(x_F)$ and $c$ coincide on $[n_0,\infty)$ and hence on $F$.
This shows that $c$ is in the closure of $\sigma'(V^\Z)$.
\par 
However, $c$ is not in the image of $\sigma'$.
Indeed, suppose on the contrary that $c = \sigma'(x)$ for some $x \in V^\Z$.
This means that $x(n + 1) = v_1+ \psi(x(n))$ for all $n \in \Z$.
By induction, we  get
$$
x(n) = v_1+ v_2 + \cdots + v_i + \psi^{i }(x(n - i ))
$$
for all $n \in \Z$ and $i \geq1$.
As the image of $\psi^{i}$ is the vector subspace of $V$ spanned by the vectors $v_{i + 1},v_{i + 2},\dots$, this implies that  $x(n) \in E$ and that the  $i$-th coordinate of $x(n)$ in  the basis of $E$ formed by the sequence $(v_i)_{i \geq 1}$ is equal to $1$.
This gives us a contradiction since the number of nonzero coordinates of a vector with respect to a given basis must always be finite.
Therefore, $c$ does not belong to the image of $\sigma'$. 
This shows that $\sigma'(V^\Z)$ is not closed in $V^\Z$.
\end{proof}

\begin{proof}[Proof of Theorem \ref{t:lca-not-closed-image}]
Suppose that $G$ is not periodic and let $H \subset G$ be an infinite cyclic subgroup.
Consider the cellular automaton $\sigma' \colon V^H \to V^H$ studied in Lemma \ref{l:sigma-2}. By virtue of this lemma and assertions (iv) and (iii) in Theorem  \ref{t:induction}, the  cellular automaton $\tau'  \colon V^G \to V^G$ induced by $\sigma'$ is linear and its image $\tau'(V^G)$ is not closed in $V^G$.
\end{proof}

% SECTION 6
\section{Concluding remarks and questions}
\label{sec:questions}

Recall that a group is called \emph{locally finite} if all its finitely generated subgroups are finite.
It is clear that every locally finite group is periodic.

\begin{proposition}
Let $G$ be a locally finite group and let $A$ be a set. Then:
\begin{enumerate}[\rm (i)]
\item
every bijective cellular automaton $\tau \colon A^G \to A^G$ is reversible;
\item
every cellular automaton $\tau \colon A^G \to A^G$ has the closed image property.
\end{enumerate}
\end{proposition}

\begin{proof}
Let $\tau \colon A^G \to A^G$ be a cellular automaton and let $M \subset G$ be a memory set for $\tau$.
Denote by $H$ the subgroup of $G$ generated by $M$ and consider the cellular automaton 
$\tau_H \colon A^H \to A^H$   over $H$ obtained by restriction of $\tau$.
\par
Since $H$ is finite, the prodiscrete topology on $A^H$ coincides with the discrete topology.
As every subset of a discrete topological space is closed, it follows that $\tau_H(A^H)$ is closed in 
$A^H$.
We deduce that $\tau(A^G)$ is closed in $A^G$ by using assertion (iii) in Theorem \ref{t:induction}.
This shows (ii).
\par
Suppose now that $\tau$ is bijective. Then $\tau_H$ is bijective by assertion (i) in Theorem \ref{t:induction}.
As $H$ is finite, every map $f \colon A^H \to A^H$ which commutes with the $H$-shift is a cellular automaton over $H$ (with memory set $H$ and local defining map $f$).
It follows that $\tau_H^{-1}$ is a cellular automaton and hence $\tau_H$ is reversible.
By applying Theorem \ref{t:induction}.(ii), we conclude that $\tau$ is reversible.
This shows (i).
\end{proof}

In view of the preceding proposition, it is natural to ask
whether the results presented in Section \ref{sec:introduction} extend to periodic groups which are not locally finite. 
 Note that, by Theorem \ref{t:induction},
it suffices to consider finitely-generated infinite periodic groups.
Examples of such groups are provided by the free Burnside groups $B(m,n)$ on $m \geq 2$ generators with large odd exponent  $n$    and the celebrated Grigor{\v{c}}uck group \cite{grig-burnside} which is a infinite $2$-group generated by $3$ involutions.

   %%% REFERENCES
\bibliographystyle{siam}
\bibliography{periodic}

\end{document}